\documentclass{article}
\usepackage{spconf,amsmath,graphicx}
\usepackage{cite}
\usepackage{amssymb}
\usepackage{amsthm}
\usepackage{graphicx}
\usepackage{mathrsfs}
\usepackage{bm}
\usepackage{caption}
\usepackage{booktabs}
\usepackage{amsmath}
\usepackage{mathrsfs}
\usepackage{mathtools}
\usepackage[dvipsnames]{xcolor}
\usepackage{float}   
\usepackage{cleveref}

\newcommand{\diag}{\mbox{diag}}

\DeclarePairedDelimiter{\floor}{\lfloor}{\rfloor}
\DeclarePairedDelimiter{\bra}{\lbrack}{\rbrack}
\DeclarePairedDelimiter{\bre}{\lbrace}{\rbrace}
\DeclarePairedDelimiter{\para}{(}{)}
\newcommand{\N}[0]{\mathbb{N}}
\newcommand{\R}[0]{\mathbb{R}}

\makeatletter
\newcommand*\bigcdot{\mathpalette\bigcdot@{0.5}}
\newcommand*\bigcdot@[2]{\mathbin{\vcenter{\hbox{\scalebox{#2}{$\m@th#1\bullet$}}}}}
\makeatother
\DeclareMathOperator{\Tr}{Tr}
\DeclareMathOperator{\E}{\mathbb{E}}
\DeclareMathOperator{\Prob}{\mathbb{P}}

\title{Information Flow Optimization in Inference Networks}
\name{Aditya Deshmukh$^{\star }$, Jing Liu$^{\star }$, Venugopal V. Veeravalli$^{\star }$, Gunjan Verma$^{\dagger}$
\address{$^{\star }$University of Illinois at Urbana-Champaign \ $^{\dagger}$Army Research Laboratory}}
\address{}
\begin{document}
%
\maketitle
\begin{abstract}
  The problem of maximizing the information flow through a sensor network tasked with an inference objective at the fusion center is considered. The sensor nodes take observations, compress and send them to the fusion center through a network of relays. The network imposes capacity constraints on the rate of transmission in each connection and flow conservation constraints. It is shown that this rate-constrained inference problem can be cast as a Network Utility Maximization problem by suitably defining the utility functions for each sensor, and can be solved using existing techniques. Two practical settings are analyzed: multi-terminal parameter estimation and binary hypothesis testing. It is verified via simulations that using the proposed formulation gives better inference performance than the Max-Flow solution that simply maximizes the total bit-rate to the fusion center.
\end{abstract}
\begin{keywords} 
Flow optimization, sensor networks, statistical inference, Internet of Things
\end{keywords}
\section{Introduction}
\label{sec:intro}

 
Consider a sensor network comprised of sensor nodes, relay/intermediate nodes and a fusion center (FC). Each sensor node takes a single observation, which is a random variable related to the state of nature/target of interest. Assume that the sensor nodes cannot communicate with each other, and that they compress their observation and send it to the intermediate nodes.  The intermediate  nodes  just  act  as  relays  and  route  messages towards  the  fusion center. The fusion center is tasked with inferring the state of nature. The relay network is considered to be capacitated and satisfies flow conservation at each relay node, i.e. number of units of information received is equal to the number of units of information transmitted by each relay node. So, the problem is a rate-constrained inference problem. Note that the trade-off here is between the rate assigned to each sensor versus the quality of compressed observations. We show that if the inference objective can be modeled as a sum of utility functions, which capture the information content in the compressed observation of each sensor according to the rate assigned to them, then the problem can be cast as a Network Utility Maximization problem. We analyze two inference tasks: parameter estimation and detection (specifically binary hypothesis testing), and explicitly formulate the utility functions for these tasks.

There is considerable literature on multi-terminal parameter estimation, either from one-shot or having time series (e.g.,~\cite{623151,669162,1365154,ekrem2014outer,4244703,4355257,amari1998statistical}). To the best of our knowledge, they at most consider the sum rate constraint for the sensors or characterizing their rate region for a given distortion constraint, without considering the actual network structure and link capacity constraints (e.g., whether the desired rates for each sensor are achievable by the actual network).

We focus on the one-shot problem in this paper. The closet work to our multi-terminal (vector) parameter estimation is Sani and Vosoughi~\cite{7450677}, in which it is assumed that the unknown vector is zero-mean Gaussian with known covariance matrix.  We do not need this assumption and treat the vector as deterministic and unknown. Another distinguishing factor is that \cite{7450677}  assumes there is an overall sum rate constraint on the network for these sensors, while our information flow optimization scheme takes into account the actual structure and link capacity constraints of the network.

Prior work that studies detection in sensor networks, includes Tay and Tsitsiklis \cite{Tay2008}, in which the one-shot problem of decentralized detection in tree structured sensor networks is considered. They characterized the optimal error exponent under a Neyman-Pearson formulation, but did not analyze capacitated networks. Another line of work was initiated by Han and Amari \cite{han}, in which the problem of multi-terminal hypothesis testing is studied with one-hop rate constraints to a fusion center and sequential observations. The  emphasis of this work is on extending rate distortion  theory to this distributed setting.

\section{Problem Setup}
\label{sec:setup}
 We model the network as an undirected graph $G=(V,E)$. In this graphical model, each edge $(u,v)\in E$ has an integral capacity $c_{uv}$ associated with it, which is the maximum possible rate in that connection(either way). Let the rate in edge $(u,v)\in E$ be $r_{uv}$. We define the rates to be anti-symmetric, i.e. $r_{uv}=-r_{vu}$. The rates obey flow conservation at each intermediate node. Let $S$ denote the set of sensor nodes. Let $N=|S|$ be the number of sensor nodes. Let $t$ denote the fusion center. The performance of the graph is measured in terms of the information content regarding the state of nature (to be inferred) in the messages sent by the sensor nodes. In this regard we assume there exists a concave utility function $g_s:\R\to\R$, that captures the notion of information content in the message from sensor node $s$. The goal is to maximize the sum utilities of the messages received at the fusion center. Let the sum rate from a node $u\in V$ be denoted by 
\begin{equation}
r_u = \sum_{v\in V: (u,v)\in E}{r_{uv}}
\end{equation}
We pose the optimization problem as follows
\begin{subequations}
\label{problem}
\begin{align}
    &\max_{r_{uv}} \sum\limits_{s\in S} g_s(r_s)\label{eq:P_obj}\\
    \text{subject to } &\forall (u,v)\in E, \;r_{uv}\leq c_{uv}\label{eq:cap}\\
                       &\forall (u,v)\in E, \;r_{uv}= -r_{vu}\label{eq:anti}\\
                       &\forall u\in V\setminus \bre{S\cup t}, \; r_u =0.\label{eq:cons}
\end{align}
\end{subequations}

In this problem formulation, we only consider the one-shot problem, i.e. only a single observation (or a single block of observations) is collected, compressed and sent through the network by each sensor. Assuming the rate is real-valued, Problem \eqref{problem} is a convex program and therefore is straightforward to solve. It is worth mentioning that, many decentralized or distributed methods have been proposed to solve this kind of Network Utility Maximization problem, e.g.,~\cite{palomar2006tutorial,4110294,1551078,1093711} and the references therein. For example, it is straightforward to use Dual Decomposition method to solve the above convex problem~\cite{palomar2006tutorial,turowskamaximization}. If the rates and capacities are constrained to be integer-valued, we can use the solution of the real-valued relaxation to obtain an integral solution in polynomial time ~\cite[Theorem 5]{lee2013new}, such that the integral rates $r_{uv}^{(I)}$ are close to the real-valued rates $r_{uv}$ in the sense that $\sum\limits_{s\in S} r_s^{(I)}=\floor*{\sum\limits_{s\in S}r_s}$.

We now discuss two instances of inference tasks in sensor networks: multi-terminal parameter estimation and binary hypothesis testing.

\section{Flow Optimization for Parameter Estimation}
\label{sec:estimation}
We consider a sensor network with $N$ spatially distributed sensors. Each sensor $i$ takes a linear measurement of the unknown deterministic vector $\bm x$, and quantizes the corresponding scalar measurement $y_i$ into certain number of bits. Then the sensors want to send those bits to the fusion center through a capacitated relay network. Finally, the fusion center estimates $\bm x$ based on the quantized data it received from all the sensors. Mathematically, we assume the following linear observation model on each sensor $i$: 
\begin{align}\label{model}
y_i= \bm a_i^T \bm x+ \eta_i, i=1,...,N, 
\end{align}
\noindent which can be jointly written as:
\begin{align}\label{model2}
\bm y= \bm A \bm x+ \bm \eta,
\end{align}

\noindent where the vector $\bm y\in \mathbb{R}^{N}$ are the overall measurements from all the $N$ sensors, and the sensing matrix $\bm A\in \mathbb{R}^{N\times q}$ is known. The vector $\bm \eta \in \mathbb{R}^{N}$ models the i.i.d. {\em bounded} noise with zero mean and variance $\sigma^2$.

For simplicity, we assume that the value of $y_i$ is within the sensor input range, and the quantizer uniformly divides that range into small intervals for quantization. For each $ y_i$, we denote its quantized version as $d_i$, and the corresponding quantization noise is $ \epsilon_i =  d_i- y_i$. So we have
\begin{align}\label{model3}
d_i= \bm a_i^T \bm x+  \eta_i +  \epsilon_i, i=1,...,N, 
\end{align}
For small quantization intervals $\Delta_i$, the quantization noise $\epsilon_i$ is approximately uncorrelated with the input $ y_i$ and has zero-mean with variance $\Delta_i^2/12$~\cite{widrow2008quantization}.

The fusion center wants to accurately estimate $\bm x$ from all the quantized data it received via the least squares $\bm \hat{x}=\bm A^{\dagger} \bm d $. The question is:
\textit{Given the limited capacity of the network, to minimize the estimation error at the fusion center, how should we optimize the network flows and how many bits should we allocate to each sensor? }

The well studied max-flow method may not give the optimal solution for this problem, even though it finds the maximal possible number of bits that can be sent from those sensors to the FC.

Our goal here is not to minimize the reconstruction error of $\bm y$, but rather minimizing the mean squared error $\E[\| \hat{\bm x}-\bm x\|_2^2]$ via optimizing the network \textit{information} flows. We have that
\begin{align}
 &\E[\| \bm \hat{\bm x}-\bm x\|_2^2]\\
 =&\E[\| \bm A^{\dagger} \bm d-\bm x\|_2^2]\\
=&\E[\| \bm A^{\dagger} (\bm {Ax+\eta+\epsilon})-\bm x\|_2^2]\\
=&\E[\| \bm A^{\dagger} (\bm {\eta+\epsilon})\|_2^2]\\
=&\E[\Tr \{\bm A^{\dagger} (\bm {\eta+\epsilon})(\bm A^{\dagger} (\bm {\eta+\epsilon}))^T \}]\\
=&\E[\Tr \{\bm A^{\dagger} (\bm {\eta+\epsilon}) (\bm {\eta+\epsilon})^T(\bm A^{\dagger})^T \}]\\
=&\Tr \{\bm A^{\dagger} \E[(\bm {\eta+\epsilon}) (\bm {\eta+\epsilon})^T](\bm A^{\dagger})^T \}\\
\approx &\Tr \{\bm A^{\dagger} \diag(\sigma^2+\Delta_1^2/12,...,\sigma^2+\Delta_N^2/12)(\bm A^{\dagger})^T \}\\
=&\sum_{i=1}^N \left [(\sigma^2+\Delta_i^2/12)\|\bm A^{\dagger}(:,i)\|_2^2\right ] \label{mse}
\end{align}

Since $\Delta_i \propto \frac{1}{2^{ r_i}}$, where $ r_i$ is the number of bits allocated to sensor $i$, minimizing (\ref{mse}) is equivalent to maximizing the following:
\begin{align}
    -\sum_{i=1}^N \|\bm A^{\dagger}(:,i)\|_2^2 / 4^{r_i}, \label{mse2}
\end{align}
\noindent Now the overall problem is cast as Problem \eqref{problem}.  

\section{Flow Optimization for Detection}
\label{sec:detection}
We first discuss some known results of optimal quantizers for detection in the classical setting before introducing the sensor network setting. Let $(\mathcal{Y}, \mathcal{F})$ be a measurable space. Let $\Prob_0$ and $\Prob_1$ be two probability measures defined on $(\mathcal{Y}, \mathcal{F})$. Under hypothesis $H_i$, $i=0,1$, a $\mathcal{Y}$-valued random variable $Y$ has the distribution $\Prob_i$. We define a deterministic $n$-level quantizer as an $\mathcal{F}$-measurable function that maps $\mathcal{Y}\to [n]\coloneqq\bre{1,2,\dots,n}$.  Let $\Gamma_n$ be the set of all randomized n-level quantizers. Let the distribution of $\gamma(Y)$ under hypothesis $H_i$ be $Q_i(\gamma)$, for $\gamma \in \Gamma_n$ for some $n$. Consider the problem of finding a randamized $n$-level quantizer to maximize the Kullback-Leibler(KL) divergence between $Q_1(\gamma)$ and $Q_0(\gamma)$, denoted by $D(Q_1(\gamma)||Q_0(\gamma))$. It is well-known that the optimal error exponent achieved in Neyman-Pearson testing with access to quantized observations is $D(Q_1(\gamma)||Q_0(\gamma))$. In this regard, we define a utility function $f_{\Prob_1,\Prob_0}:\N\to\R$ as follows:
\begin{equation}\label{eq:concave}
    f_{\Prob_1,\Prob_0}(n) \coloneqq \sup_{\gamma\in\Gamma_n}D(Q_1(\gamma)||Q_0(\gamma))
\end{equation}

Assume for simplicity that $\Prob_0$ and $\Prob_1$ have densities $p_0$ and $p_1$ with respect some common measure $\mu$. The likelihood ratio is a measurable function $L:\mathcal{Y}\to \bra{0,\infty}$ given by
\begin{equation}
    L(y) = 
    \begin{cases}
        p_1(y)/p_0(y), &\;p_1(y)\neq 0, p_0(y)\neq 0\\
        0, &\; p_1(y)= 0, p_0(y)\neq 0\\
        \infty &\;p_1(y)\neq 0, p_0(y)= 0.
    \end{cases}
\end{equation}

We define the \textit{threshold set} $T_n$ to be the set of vectors $\bm{t}=\para{t_0,t_1,\dots,t_n}\in [0,\infty]^{n+1}$ satisfying $0=t_0\leq t_1\leq \dots\leq t_{n-1}\leq t_n= \infty$.
We define a \textit{likelihood ratio quantizer} $\gamma_{\bm{t}}$ with threshold vector $\bm{t}\in T_n$, as a quantizer that satisfies the following $\forall \ell \in [n], \forall i\in\bre{0,1}$:
\begin{equation}
    \Prob_i\para*{\gamma_{\bm t}(Y)=\ell \text{ and } L(Y)\notin [t_{\ell-1},t_\ell]}=0.
\end{equation}
Tsitsiklis \cite{tsitsiklis} proved that the optimal quantizer for maximizing $D(Q_1(\gamma)||Q_0(\gamma))$ for a given number of quantization levels is the likelihood ratio quantizer, i.e.
\begin{equation}\label{eq:tsitsiklis}
    \max_{\gamma\in\Gamma_n}D(Q_1(\gamma)||Q_0(\gamma))=\max_{\bm{t}\in T_n}D(Q_1(\gamma_{\bm t})||Q_0(\gamma_{\bm t})).
\end{equation}
Thus, we get $f_{\Prob_1,\Prob_0}(n) = \mathop{\max}\limits_{\bm{t}\in T_n}D(Q_1(\gamma_{\bm t})||Q_0(\gamma_{\bm t})).$

We now consider binary hypothesis detection in the sensor network setting.  Let $\Prob_i^{(j)}$ denote the distribution of observation $Y_j\in \mathcal{Y}_j$ from sensor $j$ under the hypotheses $H_i$, $i=0,1$. As before, assume that $\Prob^{(j)}_0$ and $\Prob^{(j)}_1$ are absolutely continuous with respect to some common measure $\mu^{(j)}$. The fusion center is tasked with the objective of finding the true hypothesis. In this case, the Neyman-Pearson problem with the rate conservation and capacity constraints is an intractable problem.
We consider the problem of maximizing the KL-divergence between the distributions of the quantizer outputs under the two hypothesis,
\begin{align}\label{eq:np}
    &D\para*{Q_{1}^{(1)}\otimes Q_{1}^{(2)}\otimes\dots Q_{1}^{(n)}||Q_{0}^{(1)}\otimes Q_{0}^{(2)}\otimes\dots Q_{0}^{(N)}}\\
    &=\sum\limits_{j=1}^N D\para*{Q_{1}^{(j)}||Q_{0}^{(j)}},
\end{align}
where $Q^{(j)}_i$ is a shorthand notation for $Q^{(j)}_i(\gamma_j)$, the distribution of the quantized observation $\gamma_j(Y_j)$ from the $j^\text{th}$ sensor. Note that \eqref{eq:np} is the optimal error exponent achieved in Neyman-Pearson testing with access to quantized observations from all sensors.

\begin{figure}[H]
\centering
\includegraphics[width=8.5cm]{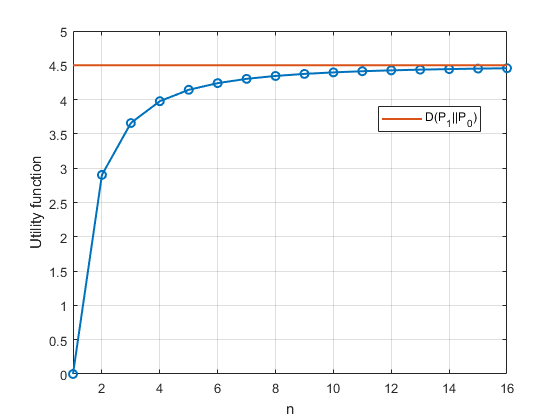}
\caption{Plot of $f_{\Prob_1,\Prob_0}(n)$, $p_0$ is $\mathcal{N}(0,1)$ and $p_1$ is $\mathcal{N}(3,1)$.}
\label{fig:1}
\end{figure}

\begin{figure}[H]
\centering
\includegraphics[width=8.5cm]{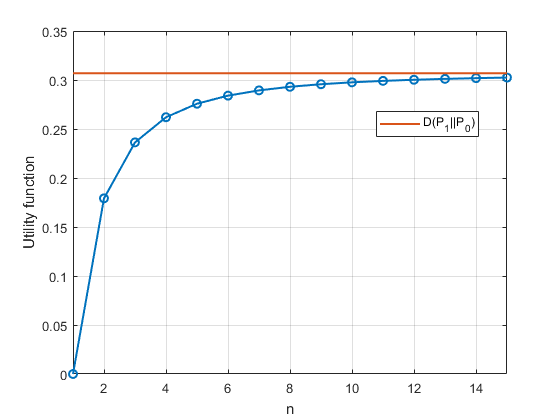}
\caption{Plot of $f_{\Prob_1,\Prob_0}(n)$, $p_0$ is Exp$(0.5)$ and $p_1$ is Exp$(1)$.}
\label{fig:2}
\end{figure}

In light of \eqref{eq:tsitsiklis}, our objective can be re-written as follows:
\begin{align}
    &\max_{r_j,\gamma_j\in\Gamma_{r_j}} \sum\limits_{j=1}^N D\para*{Q_{1}^{(j)}(\gamma_j)||Q_{0}^{(j)}(\gamma_j)}\\
    &=\max_{r_j,\bm{t_{j}}\in T_{r_j}} \sum\limits_{j=1}^N D\para*{Q_{1}^{(j)}(\gamma_{\bm{t_j}})||Q_{0}^{(j)}(\gamma_{\bm{t_j}})}\\
    &=\max_{r_j} \sum\limits_{j=1}^N f_{\Prob_1^{(j)},\Prob_0^{(j)}}(r_j)\label{eq:obj}
\end{align}
Now the problem posed can be cast as Problem \eqref{problem} with $g_s(r_s)=f_{\Prob_1^{(s)},\Prob_0^{(s)}}(r_s)$, where $r_s=\log_2(n_s)$ and $n_s$ are the number of quantization levels for sensor $s$. Note that the domain of the function \eqref{eq:concave} is discrete and so we use linear-interpolation to get a surrogate objective function for the real-valued relaxation of Problem \eqref{problem}. It can be easily shown that \eqref{eq:concave} is increasing, but it is difficult to show that it is concave. We verified via simulations that it is indeed concave for many cases of distributions. For example, Fig. \ref{fig:1} shows the plot of $f_{\Prob_1,\Prob_0}(n)$ for the case when $p_0$ and $p_1$ are unit variance Gaussian densities with means $0$ and $3$ respectively. Fig. \ref{fig:2} shows the plot of $f_{\Prob_1,\Prob_0}(n)$ for the case when $p_0$ and $p_1$ are exponential densities with rates $0.5$ and $1$ respectively.


\section{EXPERIMENTS}
\label{sec:majhead}
We generate a random graph having the following structure: $N_1$ sensor nodes form the first layer, $N_2$ intermediate nodes form the second layer, $N_3$ intermediate nodes form the third layer, and $N_4$ intermediate nodes, which form the fourth layer, are connected to the fusion center. Each node from first layer is randomly connected to $K$ nodes from the second layer. Now each connected node in the second layer is in turn randomly connected to $K$ nodes from the third layer. Similar connections are established between the third and fourth layer.

\subsection{Information Flow Optimization for Parameter Estimation}
\label{sssec:exp_estimation}
We set the network structure parameters as follows: $N_1=10$, $N_2=50$, $N_3=30$, $N_4=10$ and $K=4$. Capacities are assigned to the edges uniformly random from $\bre{1,2,\dots,15}$. We study the sensing matrix $\bm A$ considered in~\cite{RobustEstimation}, where there are some sensors close to the target and therefore have stronger signal. The rest of sensors have less proximity to the target and the signal is weak. To simulate 4 sensors with weak signal and 6 sensors with strong signal, we first generate a 10 by 3 matrix whose entries are i.i.d. drawn from the uniform distribution $U(0,1)$, then we multiply its first 4 rows by a factor $\alpha$ to obtain the sensing matrix $\bm A \in \mathbb{R}^{10\times 3}$.  The unknown vector $\bm x \in \mathbb{R}^{3}$ is randomly drawn from $\mathcal{N}(0,I)$. The bounded noise $\bm \eta_i$ is i.i.d. generated from uniform distribution $U(-0.1,0.1)$. The uniform quantizer with range [-5,5] is employed to quantize $ y_i$ into $ r_i$ number of bits.

We first generate the network and sensing matrix $\bm A$ with $\alpha=1$, and report the Mean Squared Error (MSE) of the estimated $\hat{\bm x}$ by the proposed method~\footnote{After obtaining real-valued solutions to Problem \eqref{problem}, we simply floor them to get integral value.} under 100000 Monte Carlo runs with different realizations of $\bm x$ and $\bm \eta$. Then we multiply the first 4 rows of $\bm A$ by a factor $\alpha$ and redo the simulations. We also compare with the Max-Flow solution, which maximizes the total number of bits send through the network to the Fusion Center, and allocates the corresponding bits for the sensors to quantize their measurements. Table~\ref{MSE_compare} lists the total number of allocated bits and the MSE of the estimated $\hat{\bm x}$ under one realization of the network and $\bm A$. On this example, the proposed information flow optimization scheme performs much better than the Max-Flow solution\footnote{Sometimes the benefits of the proposed method are not such significant. It depends on the network and matrix $\bm A$.}, even though it may not allocate maximum number of bits. As the value of $\alpha$ decreases, the difference between the sensors increases, and the benefits of the proposed method become more clear. 
\begin{table}[htb]
  \centering
  \captionsetup{justification=centering}
\caption{MSE of the estimated $\hat{\bm x}$ and number of bits allocated under one realization of the network and $\bm A$.}
  \label{MSE_compare}

\begin{tabular}{ |c|c|c|c| }
 \hline
MSE/bits& $\alpha$=1 & $\alpha$=0.3&$\alpha$=0.1 \\  
  \hline
  Max-Flow & 0.2673~/~74  & 1.1939~/~74& 1.5068~/~74\\ 
  \hline
 Proposed &0.0148~/~73 &0.0230~/~73& 0.0241~/~74 \\ 
 \hline
\end{tabular}
\end{table}

\subsection{Information Flow Optimization for Binary Hypothesis Testing}
\label{sssec:exp_detection}
We set the network structure parameters as follows: $N_1=10$, $N_2=4$, $N_3=3$, $N_4=2$, and $K=1$. We consider 3 different settings. For all settings, under $H_0$, the observations from all sensors follow standard normal distribution. For setting 1, under $H_1$, the observations from first five sensors follow the distribution $\mathcal{N}(11,1)$, and the rest follow $\mathcal{N}(2,1)$. For setting 2, under $H_1$, the observation from the $j^{\text{th}}$ sensor follows the distribution $\mathcal{N}(j+1,1)$. For setting 3, under $H_1$, the observations from first five sensors follow the distribution $\mathcal{N}(2,1)$, and the rest follow $\mathcal{N}(11,1)$.  Capacities are assigned to the edges uniformly at random from $\bre*{1,2,3,4,5}$. Table \ref{table:2} lists the KL-divergence between the distributions of the quantized observations under $H_1$ and $H_0$, achieved by the proposed information flow optimization scheme and Max-Flow algorithm, for one realization of a network. In the obtained realization of the network, the first five nodes have smaller overall capacity than the rest of the nodes. In this case too, we see that the proposed information flow optimization scheme performs consistently better than the max-flow solution. 

\begin{table}[htb]
  \centering
  \captionsetup{justification=centering}
\caption{KL-divergence between the distributions of the quantized observations under $H_1$ and $H_0$.}
  \label{table:2}

\begin{tabular}{ |c|c|c|c|c|c| }
 \hline
Setting\# & 1 & 2&3 \\ 
  \hline
  Max-Flow & 49.0015& 109.6984 & 135.6004 \\ 
  \hline
 Proposed &50.9881&160.1884 & 198.6157\\ 
 \hline
\end{tabular}
\end{table}
\clearpage


\bibliographystyle{IEEEbib}
\bibliography{refs}

\end{document}